\documentclass[12pt,a4paper,leqno]{amsart}

\usepackage{anysize}
\usepackage[utf8]{inputenc}
\usepackage{amsmath}
\usepackage{amssymb}
\usepackage{array}

\usepackage{hyperref}

\newtheorem{note}{Note}

\begin{document}

\title[]{Proof of Chudnovskys' series for $1/\pi$}
\author{Jesús Guillera} 
\address{University of Zaragoza, Department of mathematics, 50009 Zaragoza (Spain)}
\email{jguillera@gmail.com}
\keywords{Hypergeometric series; Ramanujan-type series for $1/\pi$; Chudnovsky's series for $1/\pi$; Elliptic modular functions; Weber modular polynomials; Modular equations}
\subjclass[2010]{33E05, 33C05, 33C20, 11F03.}

\maketitle

\begin{abstract}
We prove rational alternating Ramanujan-type series of level $1$ discovered by the brothers David and Gregory Chudnovsky, by using a method of the author. We have carried out the computations with Maple (a symbolic software for mathematics).
\end{abstract}

\section{Introduction}

In his famous paper \cite{Ramanujan} of $1914$ Ramanujan gave a list of $17$ extraordinary formulas for the number $1/\pi$, which are of the following form
\begin{equation}\label{ramanujan-series}
\sum_{n=0}^{\infty} \frac{\left(\frac12\right)_n\left(\frac1s \right)_n\left(1-\frac1s \right)_n}{(1)_n^3} (a+bn) \, z^n \, = \frac{1}{\pi}, \qquad (c)_{_0}=1, \quad (c)_{_n}=\prod_{j=1}^{n} (c+j-1),
\end{equation}
where $s \in \{2, 3, 4, 6\}$, and $z, b, a$ are algebraic numbers. Instead of using $s$ to classify them, we will use the level $\ell$ of the family (the level of the modular forms that parametrize it). It is known that
\[
\ell=4\sin^2 \frac{\pi}{s}.
\]
The only formulas of level $\ell=1 \, (s=6)$ in the list recorded by Ramanujan are
\begin{equation}
\sum_{n=0}^{\infty} \frac{\left(\frac12\right)_n\left(\frac16 \right)_n\left(\frac56 \right)_n}{(1)_n^3} (11n+1) \left( \frac{4}{125} \right)^{\! n} = \frac{5\sqrt{15}}{6 \pi}, \label{rama-example-3} 
\end{equation}
and
\begin{equation}
\sum_{n=0}^{\infty} \frac{\left(\frac12\right)_n\left(\frac16 \right)_n\left(\frac56 \right)_n}{(1)_n^3} (133n+8) \left( \frac{4}{85} \right)^{\! 3n} = \frac{85\sqrt{255}}{54 \pi}, \label{rama-example-4} 
\end{equation}
see \cite[eq. 33 and 34]{Ramanujan}. However the most interesting series in this level are the alternating ones, which were discovered by the brothers David and Gregory Chudnovsky in $1987$ \cite{Chudnovsky}. The most impressive is
\begin{equation}\label{chud-brothers} 
\sum_{n=0}^{\infty} \frac{\left(\frac12\right)_n\left(\frac16 \right)_n\left(\frac56 \right)_n}{(1)_n^3} (545140134n+13591409) \left( \frac{-1}{53360} \right)^{\! 3n} = \frac{\sqrt{640320^3}}{12 \pi}. 
\end{equation}

\begin{note} \rm
The papers \cite{Chudnovsky} and \cite{Ramanujan} are reprinted in  \cite{Pi-source-book}: A book collecting works on the number $\pi$.
\end{note}

In this paper we will prove alternating Ramanujan-type series for $1/\pi$ of level $1$ discovered by David and Gregory Chudnovsky, by using the formulas obtained by the author in \cite{Gui-meth-rama}. The fastest of all rational series for $1/\pi$ (not only for level $1$) is (\ref{chud-brothers}), which provides approximately $\log_{10} (53360^3) \simeq 14.18$ correct digits of $\pi$ per term. In \cite{Gui-rama-3A23}, we applied our method to prove the fastest series of level $3$, an alternating one discovered by Chan, Liaw and Tan \cite{Chan-Liaw-Tan}, and in \cite{Gui-rama-2P29} we proved the fastest series due to Ramanujan \cite[eq. 44]{Ramanujan}.

\section{Elliptic modular functions}

According to the papers \cite{Gui-meth-rama} and \cite{Gui-rama-3A23}, we let
\[
F_\ell(x)={}_2F_1\biggl(\begin{matrix} \frac1s, \, 1-\frac1s \\ 1 \end{matrix}\biggm| x \biggr), \qquad \ell=4 \sin^2 \frac{\pi}{s},
\]
where $s=2,\, 3, \, 4, \, 6$, and $\ell=4, \, 3, \,2 , \, 1$ is the corresponding level. The following functions $x_{\ell}(q)$ are modular functions (in a wide sense) of levels $\ell=4, \, 2, \, 3$, respectively \cite[p. 244 and p. 261]{Cooper-book}, that parametrize $x$ in such a way that $z_\ell(q)=4x_\ell(q)(1-x_\ell(q))$ are modular functions, and $F_{\ell}(x_\ell(q))$ modular forms of weight $2$:
\[
x_{4}(q) = 16 q \prod_{n=1}^{\infty} \left( \frac{1+q^{2n}}{1+q^{2n-1}} \right)^{\! 8}, \qquad
x_{2}(q) =\frac{64 q}{64 q+{\displaystyle \prod_{n=1}^{\infty} (1+q^n)^{-24}}}, 
\]
and
\[
x_{3}(q) =\frac{27 q}{27 q + {\displaystyle \prod_{n=1}^{\infty} (1+q^n+q^{2n})^{-12}}}.
\]
It is a well known theorem  that all the elliptic modular functions are algebraically related. For example, one has \cite[p. 274]{Cooper-book}:
\[
J=\frac{1728}{4x_1(1-x_1)}=\frac{64(1+3x_2)^3}{x_2(1-x_2)^2}=\frac{27(1+8x_3)^3}{x_3(1-x_3)^3}=
\frac{16(1+14x_4+x_4^2)^3}{x_4(1-x_4)^4},
\]
where $J$ is the modular invariant:
\[
J(q) =q^{-1}+744+196884\,q+21493760\,q^2+\cdots.
\]
If we let $\beta=x_{\ell}(q), \, \alpha=x_{\ell}(q^d)$,
then a modular equation of level $\ell$ and degree $1/d$ (of $\beta$ with respect to $\alpha$) or $d$ (of $\alpha$ with respect to $\beta)$, is an algebraic relation $A(\alpha, \beta)=0$. It is important to observe that in \cite{Gui-meth-rama} we showed that all we need to know, in order to prove the Ramanujan-type series for $1/\pi$, are the modular equations satisfied by $x_\ell(q)$.

\section{Weber modular equations}

Instead of using modular equations in the R. Russel form \cite{Chan-Liaw} as we did in \cite{Gui-rama-3A23} and \cite{Gui-rama-2P29}, we will use Weber modular equations for proving the Ramanujan-type series of level $1$. If we let $\beta=x_1(q)$ and $\alpha=x_1(q^d)$, and $\Phi(u,v)$ is the Weber modular polynomial \cite{Weber-polynomial} of degree $d$, then
\[ 
\alpha(1-\alpha)=\frac{432 \, u^{24}}{(u^{24}-16)^3}, \quad \beta(1-\beta)=\frac{432 \, v^{24}}{(v^{24}-16)^3}, \qquad \Phi(u,v)=0,
\]
is a modular equation of level $1$ and degree $d$. In this paper we will apply our method to prove the alternating series of level $1$ and degrees $5, 7, 11, 17$ and $41$. Our proofs of $1A5$, $1A11$, $1A17$ and $1A41$ ($\ell$Ad), where A means alternating, are completely analogues, and for proving them we will suitably modify the polynomial $\Phi(u,v)$ into another polynomial $P(u,v)$ in order to have a modular equation of the form
\[
\alpha(1-\alpha)=\frac{432 \, u^{12}}{(u^{12}-16)^3}, \quad \beta(1-\beta)=\frac{432 \, v^{12}}{(v^{12}-16)^3}, \quad P(u,v)=0,
\]
because we have observed that by doing it the computations are simpler. For proving $1A7$ we do not modify the Weber polynomial.

\section{The formulas of our method}

From a modular equation of level $\ell$ and degree $d$ ($\alpha$ respect to $\beta$), we can derive two real Ramanujan-type series for $1/\pi$:
\[
\sum_{n=0}^{\infty} \frac{\left(\frac12\right)_n\left(\frac1s \right)_n\left(1-\frac1s \right)_n}{(1)_n^3} (a+bn) \, z^n = \frac{1}{\pi}, \qquad \ell=4 \sin^2 \frac{\pi}{s},
\]
one of positive terms $z>0$ and the other one being and alternating series $z<0$. In  \cite{Gui-meth-rama} we proved that they correspond respectively to the following sets of formulas: 
\begin{equation}\label{pos}
q =e^{-\pi \sqrt{\frac{4d}{\ell}}}, \quad z=4 \alpha_0 \beta_0, \quad b=(1-2\alpha_0) \, \sqrt{\frac{4d}{\ell}}, \quad a=-2\alpha_0 \beta_0 \frac{m'_0}{\alpha'_0} \, \frac{d}{\sqrt{\ell}},
\end{equation}
and
\begin{equation}\label{alter}
q =-e^{-\pi \sqrt{\frac{4d}{\ell}-1}}, \quad z=4 \alpha_0 \beta_0, \quad b=(1-2\alpha_0) \sqrt{\frac{4d}{\ell}-1}, \quad a=-2\alpha_0 \beta_0 \frac{m'_0}{\alpha'_0} \, \frac{d}{\sqrt{\ell}},
\end{equation}
where the multiplier $m(\alpha, \beta)$ is given by the Ramanujan formula:
\begin{equation}\label{multiplier}
m^2=\frac{1}{d} \frac{\beta(1-\beta)}{\alpha(1-\alpha)} \frac{\alpha'}{\beta'},
\end{equation}
Taking logarithms in (\ref{multiplier}) and differentiating, we get
\begin{equation}\label{der-m}
\frac{m'}{\alpha'} = \frac{m}{2 \alpha'} \left( \frac{\beta'}{\beta}-\frac{\beta'}{1-\beta} -\frac{\alpha'}{\alpha} + \frac{\alpha'}{1-\alpha} + \frac{\alpha''}{\alpha'} - \frac{\beta''}{\beta'}\right),
\end{equation}
From (\ref{multiplier}) and (\ref{der-m}), we obtain the following formulas:
\begin{equation}\label{coc-m-alfa}
\frac{\beta'_0}{\alpha'_0}=\frac{1}{d m_0^2}, \qquad \frac{m'_0}{\alpha'_0} = \frac12 \left( m_0 + \frac{1}{d m_0} \right) \frac{\alpha_0-\beta_0}{\alpha_0 \beta_0} + \frac{m_0}{2\alpha'_0}\left( \frac{\alpha''_0}{\alpha'_0} - \frac{\beta''_0}{\beta'_0}\right).
\end{equation}
Hence for proving a Ramanujan-type series for $1/\pi$ of degree $d$ one only needs to know a modular equation of that degree. When we apply our method we begin making $\beta=1-\alpha$ in the modular equation and choose a solution $\alpha_0$. If with that solution we get $|m_0| \neq 1/\sqrt{d}$ then it is not of degree $d$ and we have to try another solution. A good test to select the correct solution $\alpha_0$ was explained in \cite{Gui-meth-rama} and used in \cite{Gui-rama-3A23} and  \cite{Gui-rama-2P29}.

\section{Proofs of Chudnovskys' series for $1/\pi$}
First we will prove the alternating series of degree $17$. The proofs of the alternating series of degrees $d=5,11,41$ are completely similar. Finally we will prove the series $1A7$. Other proofs of Chudnovskys' series are given in \cite{Zhao} and \cite{Milla}. 
\subsection{Proof of the formula $1A17$}\label{A17}
We see in the tables of \cite{Gui-meth-rama} that the alternating Ramanujan-type series for $1/\pi$ of level $1$ and degree $17$ is
\begin{equation}\label{-440}
\sum_{n=0}^{\infty} \frac{\left(\frac12\right)_n\left(\frac16 \right)_n\left(\frac56 \right)_n}{(1)_n^3} (261702n+10177) \left( \frac{-1}{440} \right)^{\! 3n} = \frac{3 \cdot 440^2}{\sqrt{330} \, \pi}.
\end{equation}
With other methods one needs to use a modular equation of degree $4d-1=67$ to prove it. Here we will show how to prove it from the Weber modular equation of degree $17$:
\begin{equation}\label{weber-modeq-17}
\alpha(1-\alpha)=\frac{432 \, u^{24}}{(u^{24}-16)^3}, \quad \beta(1-\beta)=\frac{432 \, v^{24}}{(v^{24}-16)^3}, \quad \Phi_{17}(u,v)=0,
\end{equation}
see $\Phi_{17}(u,v)$ at \cite[second web-page]{Weber-polynomial}. However, we prefer to transform it in the following way: First write $\Phi_{17}(u,v)=0$ as 
\[
Q(u,v)=uv R(u,v), 
\]
where
\begin{align*}
Q(u,v) &=(u^{18}+v^{18})+17(u^{16}v^{10}+u^{10}v^{16})+119(u^{12}v^{6}+u^6v^{12})+272(u^8v^2+u^2v^8), \\
R(u,v) &=-u^{16}v^{16}-34(u^{14}v^2+u^2v^{14})+34u^{12}v^{12}+340u^{8}v^{8}+544u^4v^4-256.
\end{align*}
Squaring we have $Q^2(u,v)=u^2 v^2 R^2(u,v)$. Finally, replacing $u$ with $\sqrt{u}$ and $v$ with $\sqrt{v}$ we obtain $Q^2(\sqrt{u}, \sqrt{v})-u v R^2(\sqrt{u}, \sqrt{v})=0$. It is clear that the left hand side is a polynomial $P(u,v)$, and we obtain the following modular equation:
\begin{equation}\label{weber-modeq-transf-17}
\alpha(1-\alpha)=\frac{432 \, u^{12}}{(u^{12}-16)^3}, \quad \beta(1-\beta)=\frac{432 \, v^{12}}{(v^{12}-16)^3}, \quad P(u,v)=0.
\end{equation}
We will use (\ref{weber-modeq-transf-17}) instead of (\ref{weber-modeq-17}) because the calculations are simpler, and we will use Maple (a symbolic software for mathematics) to make those computations.
\par Our method begins taking $\beta=1-\alpha$. Hence, we have to find a solution of the system
\begin{equation}\label{system}
\frac{u^{12}}{(u^{12}-16)^3}=\frac{v^{12}}{(v^{12}-16)^3}, \qquad P(u,v)=0.
\end{equation}
We choose the following solution (check!):
\begin{align*}
u_0 &=\left( \frac{91}{1200} - \frac{3\sqrt{201}}{400}  \right)H^2+ \frac13 H -\frac23, \\
v_0 &= \left\{ \left( \frac{3\sqrt{201}}{800}-\frac{91}{2400} \right)H^2 - \frac16 H - \frac23 \right\}+ \left\{
\left( \frac{-9 \sqrt{67}}{800} + \frac{91\sqrt{3}}{2400} \right)H^2-\frac{\sqrt{3}}{6}H \right\} \, i,
\end{align*}
where
\[
H=\left( 91+9\sqrt{201} \right)^{\frac13}.
\]
Substituting in (\ref{weber-modeq-transf-17}), we get
\[
\alpha_0=\frac12-\frac{651}{193600}\sqrt{22110}, \qquad \beta_0=\frac12+\frac{651}{193600}\sqrt{22110}.
\]
Hence
\[
z_0=4\alpha_0 \beta_0 = \frac{-1}{440^3}.
\]
Then, from the formula for $b$ in (\ref{alter}), we get
\[
b=\frac{43617}{96800}\sqrt{330}.
\]
We choose $u$ as the independent variable. Differentiating $P(u,v)=0$ with respect to $u$ at $u=u_0$ we obtain $v'_0$, and differentiating twice $P(u,v)$ with respect to $u$ at $u=u_0$, we get $v''_0$.  Differentiating
\begin{equation}\label{al-be-u-v}
\alpha(1-\alpha)=\frac{432 \, u^{12}}{(u^{12}-16)^3}, \quad \beta(1-\beta)=\frac{432 \, v^{12}}{(v^{12}-16)^3},
\end{equation}
with respect to $u$ at $u=u_0$, we obtain $\alpha'_0$ and $\beta'_0$, and we get
\[
m_0=\sqrt{\frac{1}{d} \, \frac{\alpha'_0}{\beta'_0}}=\frac{\sqrt{67}}{34}+\frac{1}{34} \, i, \quad |m_0|=\frac{1}{\sqrt{17}}.
\]
Then, differentiating (\ref{al-be-u-v}) twice, we obtain $\alpha''_0$ and $\beta''_0$. Finally, from the formula for $a$ in (\ref{alter}) and the formulas (\ref{multiplier}) and (\ref{coc-m-alfa}), we obtain
\[
a=\frac{10177}{580800}\sqrt{330}.
\]

\subsection{Proof of the formulas $1A5$, $1A11$, $1A41$ and $1A7$}
The proofs of the formulas $1A5$, $1A11$, $1A41$ are completely similar to the proof of $1A17$: Modify the Weber polynomial $\Phi_d(u,v)$ in the same way that we have done in the case of degree $d=17$. For proving $1A7$ do not modified the Weber polynomial. Then, continue choosing the values of $u_0$ and $v_0$ that we indicate below.
\subsubsection{Proof of the formula $1A5$} Choose
\begin{align*}
u_0 &=\left( \frac{-1}{192} + \frac{\sqrt{57}}{64}  \right)H^2 - \frac13 H +\frac23, \\
v_0 &= \left\{ \left( \frac{-\sqrt{57}}{128} + \frac{1}{384} \right)H^2 + \frac16 H + \frac23 \right\} - \left\{
\left( \frac{\sqrt{3}}{384} - \frac{\sqrt{171}}{128} \right)H^2-\frac{\sqrt{3}}{6}H \right\} \, i,
\end{align*}
where
\[
H=\left(1+3\sqrt{57} \right)^{\frac13}.
\]
Then follow the steps of (\ref{A17}).
\subsubsection{Proof of the formula $1A11$} Choose
\begin{align*}
u_0 &=\left( \frac{-35}{48} + \frac{\sqrt{129}}{16}  \right)H^2 - \frac13 H + \frac43, \\
v_0 &= \left\{ \left( \frac{35}{96}-\frac{\sqrt{129}}{32} \right)H^2 + \frac16 H + \frac43 \right\} + \left\{
\left( \frac{-35 \sqrt{3}}{96} + \frac{3\sqrt{43}}{32} \right)H^2+\frac{\sqrt{3}}{6}H \right\} \, i,
\end{align*}
where
\[
H=\left(35+3\sqrt{129} \right)^{\frac13}.
\]
Then follow the steps of (\ref{A17}).
\subsubsection{Proof of the formula $1A41$} Choose
\begin{align*}
u_0 &=\left( \frac{-467}{13872} + \frac{11\sqrt{489}}{4624}  \right)H^2 - \frac13 H + \frac43, \\
v_0 &= \left\{\left(\frac{467}{27744}-\frac{11\sqrt{489}}{9248} \right)H^2 + \frac16 H + \frac43 \right\} + \left\{
\left( \frac{467 \sqrt{3}}{27744} - \frac{33\sqrt{163}}{9248} \right)H^2-\frac{\sqrt{3}}{6}H \right\} \, i,
\end{align*}
where
\[
H=\left(467+33\sqrt{489} \right)^{\frac13}.
\]
Then follow the steps of (\ref{A17}). 

\begin{note} \rm
A Maple program which automatically proves the series $1A5$, $1A11$, $1A17$ and $1A41$ is available at the web-site of the author. The procedure $\texttt{chud($\cdot$)}$ does it. For example $\texttt{chud(41);}$ automatically proves the Chudnovskys' series (\ref{chud-brothers})
\end{note}

\subsubsection{Proof of the formula $1A7$} 
Take $P(u,v)=\Phi_7(u,v)$ (that is, do not modify the Weber polynomial), and choose
\begin{align*}
u_0 &= \sqrt[6]{2} \left(\sqrt[3]{2}-1\right)^{\frac13}, \\
v_0 &= \sqrt[6]{2} \left(\sqrt[3]{2}-1\right)^{\frac13} \left( \frac12 (\sqrt[3]{4}+1) - \frac{\sqrt{3}}{6} 
(1+\sqrt[3]{4}+2\sqrt[3]{2}) \, i \right).
\end{align*}
Then follow the steps of (\ref{A17}).

\end{document}